\documentclass[11pt]{article}
\usepackage{amssymb}
\usepackage{color}

\normalbaselines

\newtheorem{Theorem}{Theorem}[section]

\newtheorem{Proposition}[Theorem]{Proposition}

\newtheorem{Lemma}[Theorem]{Lemma}

\newtheorem{Corollary}[Theorem]{Corollary}

\newtheorem{Remark}[Theorem]{Remark}

\newcommand{\RR}{{{\rm I} \kern -.15em {\rm R} }}

\newcommand{\C}{{{\rm l} \kern -.42em {\rm C} }}

\newcommand{\nat}{{{\rm I} \kern -.15em {\rm N} }}

\newcommand{\be}{\begin{equation}}
\newcommand{\ee}{\end{equation}}
\newcommand{\beq}{\begin{eqnarray}}
\newcommand{\eeq}{\end{eqnarray}}
\newcommand{\beqs}{\begin{eqnarray*}}
\newcommand{\eeqs}{\end{eqnarray*}}
\newcommand{\bt}{\begin{Theorem}}
\newcommand{\et}{\end{Theorem}}
\newcommand{\br}{\begin{Remark}}
\newcommand{\er}{\end{Remark}}
\newcommand{\bc}{\begin{Corollary}}
\newcommand{\ec}{\end{Corollary}}
\newcommand{\bl}{\begin{Lemma}}
\newcommand{\el}{\end{Lemma}}
\newcommand{\bd}{\begin{definition}}
\newcommand{\ed}{\end{definition}}

\title{Stability results for abstract evolution equations with\\ intermittent time-delay feedback}
\author{
{\sc Cristina Pignotti}
\\Dipartimento di Ingegneria e Scienze dell'Informazione e Matematica\\
 Universit\`{a} di L'Aquila\\
Via Vetoio, Loc. Coppito, 67010 L'Aquila Italy}

\date{}

\begin{document}

\textwidth=160 mm

\textheight=225mm

\parindent=8mm

\frenchspacing

\maketitle

\begin{abstract}
We consider abstract evolution equations with on--off time delay feedback. Without the time delay term, the model is described by an exponentially stable semigroup. We show that, under appropriate conditions involving the delay term, the system remains asymptotically stable. Under additional assumptions exponential stability results are also obtained. Concrete examples illustrating the abstract results are finally given.
\end{abstract}

\vspace{5 mm}

\def\qed{\hbox{\hskip 6pt\vrule width6pt
height7pt
depth1pt  \hskip1pt}\bigskip}



\section{Introduction}
\label{pbform}\hspace{5mm}

\setcounter{equation}{0}

In this paper we study the stability properties of abstract evolution equations in presence of a time delay term.

In particular, we include into the model an on--off time delay feedback, i.e. the time delay is intermittently present.

Let ${\mathcal H}$ be a Hilbert space, with norm $\Vert\cdot\Vert ,$ and let ${\mathcal A}: {\mathcal H} \rightarrow {\mathcal H}$  be a dissipative  operator
generating a $C_0-$ semigroup  $(S(t))_{t\ge 0}$ exponentially stable,
namely there are two positive constants $M$ and 
$\mu$ such that 
\begin{equation}\label{semigroup}
\Vert S(t)\Vert_{{\mathcal L}({\mathcal H})}\le Me^{-\mu t}\,,\quad \forall\ t\ge 0\,,
\end{equation}
where ${{\mathcal L}({\mathcal H})}$
denotes the space of bounded linear operators from ${\mathcal H}$ into itself.

We consider the following problem
\begin{equation}\label{problema}
\left\{
\begin{array}{l}
U_t(t)={\mathcal A}U(t)+{\mathcal B}(t) U(t-\tau ) \quad t >0,\\
U(0)=U_0,
\end{array}
\right.
\end{equation}
where $\tau\,,$ the time delay, is a fixed positive constant, the initial datum $U_0$ belongs to ${\mathcal H}$  and, for $t>0,$ ${\mathcal B}(t)$ is a bounded operator from ${\mathcal H}$ to   ${\mathcal H}.$ 

In particular, we assume 
that there exists an increasing sequence of positive real numbers $\{t_n\}_n,$ with $t_0=0,$ 
such that
$$
\begin{array}{l}
{ 1)} \quad {\mathcal B}(t)=0\quad  \forall\ t\in I_{2n}=[t_{2n},t_{2n+1}),\\
{2)} \quad  {\mathcal B}(t)={\mathcal B}_{2n+1} 
\quad \forall \  t\in I_{2n+1}=[t_{2n+1},t_{2n+2}).
\end{array}
$$
We denote $B_{2n+1}=\Vert {\mathcal B}_{2n+1} 
\Vert_{{\mathcal L}({\mathcal H})}\,,$ $n\in\nat\,.$
Moreover, denoted by $T_n$ the length of the interval $I_n,$ that is
\begin{equation}\label{Tn}
T_n=t_{n+1}-t_n,\quad n\in \nat\,,
\end{equation}
we  assume 
\begin{equation}\label{T2n}
 T_{2n}\geq \tau , \quad  \forall\  n\in\nat\,.
\end{equation}

Time  delay
effects
are frequently present
in applications and concrete models and it is now well--understood that even an arbitrarily small delay in the feedback may destabilize a system which is uniformly  stable in absence of delay (see e.g. \cite{Datko, DLP, NPSicon06, XYL}).

We want to show that, under appropriate assumptions involving  the delay feedback coefficients, the size of the time intervals where the delay appears and  the parameters $M$ and $\mu$ in (\ref{semigroup}),   the considered model is asymptotically stable or exponentially stable, in spite of the presence of the time delay term.

Stability results for second--order evolution equations with intermittent damping
were first studied by Haraux, Martinez and Vancostenoble \cite{HMV}, without any time delay term. They considered a model with intermittent on--off or with positive--negative damping
and gave sufficient conditions ensuring that the  behavior of the system in the time intervals with the standard dissipative damping, i.e.  with positive coefficient,  prevails  over the {\em bad} behavior in remaining intervals where
the damping is no present or it is present with the negative sign, namely as anti--damping.
Therefore, asymptotic/exponential stability results were obtained. 

More recently Nicaise and Pignotti \cite{ADE2012, JDDE14} considered
second--order evolution equations with intermittent delay feedback. These results have been improved and extended to some semilinear equations in \cite{GenniCri}.
In the studied models, when the delay term (which possess a destabilizing effect) is not present,
a not--delayed damping acts. Under  appropriate sufficient  conditions, stability results are then obtained.
Related results for wave equations with intermittent delay feedback have been obtained, in 1-dimension, in \cite{Gugat}, \cite{GT} and \cite{ANP2012} by using a different approach based on  the D'Alembert formula. However, this last approach furnishes stability results  only for particular choices of the time delay.

In the recent paper \cite{pignotti2016}, the intermittent delay feedback is compensated by a viscoelastic damping with exponentially decaying kernel.

The asymptotic behavior of  wave--type equations  with infinite memory and time delay feedback  has been studied
by Guesmia in \cite{Guesmia} via a Lyapunov approach and by Alabau--Boussouira, Nicaise and Pignotti \cite{AlNP2015} by combining multiplier identities and perturbative arguments.

We refer also to Day and Yang
 \cite{DY} for the same kind of problem  in the case of finite memory.
 In these papers the authors prove exponential stability results if the coefficient of the delay damping is sufficiently small. These stability results could be easily extended to
a variable coefficient $b(\cdot )\in L^\infty (0, +\infty )$ under a suitable {\em smallness} assumption on the  $L^\infty-$ norm of $b(\cdot )\,.$

In \cite{pignotti2016}, instead,
asymptotic stability  results are obtained without smallness conditions related to the $L^\infty-$ norm
of the delay coefficient.
On the other hand, the analysis is restricted to intermittent delay feedback. 
Asymptotic stability is proved when the coefficient of the delay feedback belongs to $L^1(0, +\infty )$ and the length of the time intervals where the delay is not present is sufficiently large.  The same paper considers also problems with on--off anti--damping instead of a time delay feedback. Stability results are obtained even in this case under analogous assumptions.

The idea is here to generalize the results of \cite{pignotti2016} by considering abstract evolution equations for which, without considering  the intermittent delay term, the associated operator generates an exponentially stable $C_0-$ semigroup.  

For such a class of evolution equations we already know that, under a suitable smalness condition on the delay feedback coefficient, an exponential stability result holds true (see \cite{JEE15}).
We want to show that stability results are avalaible also under a condition on the $L^1-$norm of the delay coefficient, without restriction on the pointwise $L^\infty-$norm.

The paper is organized as follows. In section \ref{well}   we give a well--posedness result. In sections \ref{st} and \ref{stexp} we prove
asymptotic and exponential stability results, respectively,  for the abstract model
 under appropriate conditions.  Stability results
 are established also for a problem with intermittent anti--damping instead
of delay feedback in section \ref{anti}.
Finally, in section  \ref{esempi}, we give some concrete applications of the abtract results.

\section{Well-posedness \label{well}}

\hspace{5mm}

\setcounter{equation}{0}

In this section we illustrate  a well-posedness results for problem
(\ref{problema}).

\begin{Theorem}\label{texistence}
For any initial datum $U_0\in {\mathcal H}$ there exists a unique (mild) solution 
 $U\in C([0,\infty); {\mathcal H})$ of problem $(\ref{problema}).$
Moreover,
\begin{equation}\label{rappresentazione}
U(t)=S(t) U_0+\int_0^t S(t-s){\mathcal B}(s) U(s-\tau )\, ds\,.
\end{equation}
\end{Theorem}
\noindent {\bf Proof.} We prove the existence and uniqueness result on the interval $[0,t_2];$ then the global result follows by translation (cfr. \cite{ADE2012}). 
In the time interval $[0, t_1],$ since ${\mathcal B} (t)= 0$ $\forall\ t\in [0, t_1),$ then there exists a unique solution $U\in C([0, \tau ], {\mathcal H})$ satisfying $(\ref{rappresentazione}).$
The situation is different in the time interval $[t_1, t_2]$ where the delay feedback is present. In this case, we decompose the interval $[t_1, t_2]$ into the successive intervals $[t_1+j\tau, t_1+(j+1)\tau ),$ for $j=0, \dots, N-1,$ where $N$ is such that $t_1+(N+1)\tau\ge t_2\,.$
The last interval is then $[t_1+N\tau, t_2]\,.$
Now, first we look at the problem on the interval $[t_1, t_1+\tau]\,.$ Here $U(t-\tau)$ can be considered as a known function. Indeed, for $t\in [t_1, t_1+\tau]\,,$ then $t-\tau \in [0, t_1],$ and we know the solution $U$ on $[0, t_1]$ by the first step. Thus, problem $(\ref{problema})$
may be reformulated on $[t_1, t_1+\tau ]$ as
\begin{equation}\label{problema2}
\left\{
\begin{array}{l}
U_t(t)={\mathcal A}U(t)+g_0(t)
 \quad t \in (\tau , 2\tau ),\\
U(\tau )=U (\tau_-), 
\end{array}
\right.
\end{equation}
where $g_0(t)={\mathcal B}(t)U(t-\tau )\,.$ 
This problem has a  unique solution $U\in C([\tau , 2\tau ], {\mathcal H} )$ (see e.g. Th. 1.2, Ch. 6 of 
\cite{pazy}) given by
$$U(t)=S(t-\tau )U(\tau_-)+\int_{\tau }^t S(t-s)g_0(s) \, ds,\quad \forall \ t\in [\tau , 2\tau ]\,.$$
Proceedings analogously in the successive time intervals $[t_1+j\tau, t_1+(j+1)\tau ),$ we obtain a solution on $[0, t_2]\,.$\qed

\section{Asymptotic stability results
\label{st}}

\hspace{5mm}

\setcounter{equation}{0}

Let $T^*$ be defined
as

\begin{equation}\label{Tstar}
T^*:=\frac 1 {\mu } \ln M\,,
\end{equation}
where $M$ and $\mu$ are the constants in $(\ref{semigroup}),$
that is $T^*$ is the time for which $Me^{-\mu T^*}=1\,.$

We can state a first estimate on the intervals $I_{2n}$ where the delay feedback is not present.

\begin{Proposition}\label{obs}
Assume $T_{2n}>T^*\,.$ Then, there exists a constant $c_n\in (0,1)$ such that
\begin{equation}\label{stimaobs}
\Vert U(t_{2n+1})\Vert^2 \le c_n \Vert U(t_{2n})\Vert^2\,,
\end{equation}
for every solution of problem $(\ref{problema}).$
\end{Proposition}

\noindent {\bf Proof.}
Observe that in the time interval $I_{2n}=[t_{2n}, t_{2n+1}]$ the delay feedback is not present since ${\mathcal B}(t)\equiv 0\,.$
Thus, (\ref{stimaobs}) easily follows from (\ref{semigroup}) with
$\sqrt{c_n}= M e^{-\mu T_{2n}}<M e^{-\mu T^{*}}=1\,.$
\qed

Let us now introduce the Lyapunov functional
\begin{equation}\label{energy}
F(t)=F(U,t):= \frac 12\Vert U(t)\Vert^2 +\frac 12 \int_{t-\tau }^t\Vert {\mathcal B}(s+\tau )\Vert_{{\mathcal L}({\mathcal H})} \Vert U(s)\Vert^2\, ds\,.
\end{equation}

\begin{Proposition}\label{intervallicattivi}
Assume $1), 2)\,.$ Moreover,
assume $T_{2n}\ge\tau,\ \forall\ n\in \nat .$ Then,
\begin{equation}\label{stimacattiva}
F^\prime (t) \le  {B}_{2n+1} \Vert U(t)\Vert^2\,,\quad t\in I_{2n+1}= [t_{2n+1}, t_{2n+2}], \ \forall\ n\in \nat\,.
\end{equation}
for any solution of problem $(\ref{problema}).$
\end{Proposition}

\noindent {\bf Proof.}
By differentiating the energy $F(\cdot ),$ we have
$$\begin{array}{l}
\displaystyle{
F^\prime (t)=\langle U(t), {\mathcal A} U(t)\rangle + \langle U(t), {\mathcal B(t)} U(t-\tau )\rangle +\frac 1 2\Vert {\mathcal B}(t+\tau )\Vert_{{\mathcal L}({\mathcal H})} \Vert U(t)\Vert^2}\\
\hspace{5 cm}\displaystyle{
-\frac 12\Vert {\mathcal B}(t)\Vert_{{\mathcal L}({\mathcal H})} \Vert U(t-\tau )\Vert^2}\,.\\
\end{array}
$$

Then, since the operator ${\mathcal A}$ is dissipative, one can estimate 

\begin{equation}\label{cri10}
\begin{array}{l}
\displaystyle{
F^\prime (t)\le  \Vert {\mathcal B}(t)\Vert_{{\mathcal L}({\mathcal H})}
 \Vert U(t)\Vert \Vert U(t-\tau )\Vert +
\frac 1 2\Vert {\mathcal B}(t+\tau )\Vert_{{\mathcal L}({\mathcal H})} \Vert U(t)\Vert^2 }\\
\hspace{4 cm}\displaystyle{
-\frac 12\Vert {\mathcal B}(t)\Vert_{{\mathcal L}({\mathcal H})} \Vert U(t-\tau )\Vert^2\,.}
\end{array}
\end{equation}

Therefore, from Cauchy--Schwarz inequality,
$$F^\prime (s)\le 
\frac 1 2\Vert {\mathcal B}(t )\Vert_{{\mathcal L}({\mathcal H})} \Vert U(t)\Vert^2 +\frac 1 2\Vert {\mathcal B}(t+\tau )\Vert_{{\mathcal L}({\mathcal H})} \Vert U(t)\Vert^2\,.$$

Now, observe that, since $T_{2n}\ge\tau,$ for every $n\in\nat\,,$ if $t$ belongs to $I_{2n+1}$ then
$t+\tau$ belongs to $I_{2n+1}$ or to $I_{2n+2}\,.$
In the first case $\Vert {\mathcal B}(t)\Vert_{{\mathcal L}({\mathcal H})}={B}_{2n+1}$ while, in the second case 
$\Vert {\mathcal B}(t)\Vert_{{\mathcal L}({\mathcal H})}=0\,.$
Thus (\ref{stimacattiva}) is proved.\qed

\begin{Theorem}\label{CP1}
Assume $\mbox{\rm 1), 2)}$ and $T_{2n}\ge \tau$ for all $n\in\nat.$ Moreover
assume $T_{2n}> T^*,$ for all $n\in\nat,$ where $T^*$ is the time defined in $(\ref{Tstar}).$ Then, if
\begin{equation}\label{general}
\sum_{n=0}^\infty
\ln \left [
e^{2{B}_{2n+1}T_{2n+1}}
 (c_n+T_{2n+1} {B}_{2n+1} )\right ] = -\infty\,,
\end{equation}
the
equation $(\ref{problema})$ is asymptotically stable, namely any solution  $U$ of $(\ref{problema})$ satisfies
$\Vert U(t)\Vert\rightarrow 0$ for $t\rightarrow +\infty\,.$
\end{Theorem}

\noindent {\bf Proof.}
Note that from (\ref{stimacattiva}) we obtain 
$$F^{\prime}(t)\le 2{B}_{2n+1} F(t),\quad t\in I_{2n+1}=[t_{2n+1},t_{2n+2}),\ n\in\nat.$$
Then, by integrating on the time interval $I_{2n+1},$
\begin{equation}\label{cri0}
F(t_{2n+2})\le e^{2{B}_{2n+1}T_{2n+1}}F(t_{2n+1}),
\quad \forall \ n\in\nat.
\end{equation}
From the definition of the Lyapunov functional $F$,
\begin{equation}\label{cri1}
F(t_{2n+1})=\frac 1 2 \Vert U(t_{2n+1})\Vert^2+\frac 1 2 \int_{t_{2n+1}-\tau }^{t_{2n+1}}\Vert {\mathcal B}(s+\tau )\Vert_{{\mathcal L}(H)}\, \Vert  U(s)\Vert^2 ds\,.
\end{equation}
Note that, for $t\in [t_{2n+1}-\tau , t_{2n+1}),$ then $t+\tau \in [t_{2n+1}, t_{2n+1}+\tau )$ and therefore, since $\vert I_{2n+2}\vert\ge \tau\,$ it results
$t+\tau \in I_{2n+1}\cup I_{2n+2}\,.$
Now, if $t+\tau \in I_{2n+2},$  then ${\mathcal B}(t+\tau )=0.$ Otherwise, if $t+\tau \in I_{2n+1},$  then $\Vert {\mathcal B}(t+\tau)\Vert = { B}_{2n+1}\,.$
Then, from (\ref{cri1}) we deduce
\begin{equation}\label{cri2}
F(t_{2n+1})=\frac 12 \Vert U(t_{2n+1})\Vert^2+\frac 1 2  {B}_{2n+1} \int_{t_{2n+1}-\tau }^{\min ({t_{2n+2}-\tau }, t_{2n+1})}\Vert U(s)\Vert^2 ds\,,
\end{equation}
since if  $t_{2n+1}> t_{2n+2}-\tau ,$ then  ${\mathcal B}(t)=0$ for all $t\in [t_{2n+2}, t_{2n+1}+\tau )\subset [t_{2n+2}, t_{2n+3}).$

Then, since $\Vert U(\cdot )\Vert$ is decreasing in the intervals $I_{2n}$ (the operator ${\mathcal A}$ is dissipative and ${\mathcal B(t)}\equiv 0$), we deduce
\begin{equation}\label{cri3}
\begin{array}{l}
\displaystyle{
F(t_{2n+1})\le \frac 12\Vert U(t_{2n+1})\Vert^2+\frac 12 T_{2n+1} {B}_{2n+1} \Vert U(t_{2n+1}-\tau )\Vert^2}\\
\hspace{2 cm}\displaystyle{
\le \frac 12\Vert U(t_{2n+1})\Vert^2+\frac 12 T_{2n+1} { B}_{2n+1} \Vert U(t_{2n})\Vert^2
\,.}
\end{array}
\end{equation}

Using this last estimate in (\ref{cri0}), we obtain

\begin{equation}\label{cri4}
\Vert U(t_{2n+2})\Vert^2 \le 2 F(t_{2n+2})\le e^{2{B}_{2n+1}T_{2n+1}}
 (c_n+ T_{2n+1} { B}_{2n+1}
)\Vert U(t_{2n})\Vert^2,
\quad \forall \ n\in\nat,
\end{equation}
where we have used also the estimate (\ref{stimaobs}).
By iterating  this argument we arrive at

\begin{equation}\label{cri5}
\Vert U(t_{2n+2})\Vert^2\le\displaystyle{\Pi_{k=0}^n }e^{2{B}_{2k+1}T_{2k+1}}
 (c_k+T_{2k+1} {B}_{2k+1}
)\Vert U_0\Vert^2 ,
\quad \forall \ n\in\nat\,.
\end{equation}

Now observe that $\Vert U(t)\Vert $ is not decreasing in the whole $(0,+\infty).$ However, it is decreasing
for $t\in [t_{2n}, t_{2n+1})$, $n\in\nat\,,$ where the destabilizing  delay feedback does not act and so
\begin{equation}\label{Nice1}
\Vert U(t)\Vert \le \Vert U(t_{2n})\Vert ,\quad \forall \; t\in [t_{2n}, t_{2n+1}).
\end{equation}
Moreover, from (\ref{cri3}), for $t\in [t_{2n+1}, t_{2n+2})$ we have 
\begin{equation}\label{Nice2}
\Vert U(t)\Vert^2\le 2F(t)\le e^{2{B}_{2n+1}T_{2n+1}}(c_n+B_{2n+1}T_{2n+1}) \Vert U(t_{2n})\Vert^2,
\end{equation}
where in the second inequality we have used  (\ref{stimaobs}).

Then, we have asymptotic stability if

$$\displaystyle{
\Pi_{k=0}^n} e^{2{ B}_{2k+1}T_{2k+1}}
 (c_k+T_{2k+1} {B}_{2k+1}) \longrightarrow 0,\quad \mbox{\rm for}\ n\rightarrow\infty,$$
 or equivalently
 $$\ln \Big [
\Pi_{k=0}^n e^{2{ B}_{2k+1}T_{2k+1}}
 (c_k+T_{2k+1} { B}_{2k+1} )\Big ] \longrightarrow -\infty,\quad \mbox{\rm for}\ n\rightarrow\infty\,,$$
 namely under the assumption (\ref{general}).
 This concludes the proof. \qed

\begin{Remark}\label{particular}
{\rm
In particular,
 (\ref{general}) is verified if the following conditions are satisfied:
\begin{equation}\label{star1bis}
\sum_{n=0}^\infty { B}_{2n+1}T_{2n+1}<+\infty\quad \mbox{\rm and}\quad \sum_{n=0}^\infty \ln c_n=-\infty\,.
\end{equation}
Indeed, it is easy to see that (\ref{star1bis}) is equivalent to 
\begin{equation}\label{star1}
\sum_{n=0}^\infty B_{2n+1}T_{2n+1}<+\infty\quad \mbox{\rm and}\quad \sum_{n=0}^\infty \ln (c_n+B_{2n+1}T_{2n+1})=-\infty
\end{equation} and that ({\ref{star1})
implies (\ref{general}).}

Therefore, from $(\ref{star1bis}),$  we have stability if $\Vert {\mathcal B}(t)\Vert \in L^1(0,+\infty )$ and,
for instance,
the length of the {\em good} intervals $I_{2n}$ is greater than a fixed time $\bar T,$ $\bar T>T^*$ and $\bar T\ge\tau ,$
namely
$$T_{2n}\ge \bar T,\quad\forall\ n\in \nat\,.$$ Indeed, in this case there exists $\bar c\in (0,1)$ such that $0<c_n<\bar c\,.$ }
\end{Remark}

If we assume that the length
of the delay intervals, namely the time intervals where the delay feedback is present, is lower
than the time
delay $\tau,$ that is
\begin{equation}\label{rest}
T_{2n+1}\le\tau,\quad \forall n\in\nat\,.
\end{equation}
we can prove another asymptotic stability result which is, in some sense, complementary to the previous one.

In this case we can directly work with $\Vert U(t)\Vert$ instead of passing trough the  function $F(\cdot).$
We can give the following preliminary estimates
on the time intervals $I_{2n+1},$ $n\in\nat.$

\begin{Proposition}\label{CP2}
Assume $\mbox{\rm 1),\ 2)}.$ Moreover assume $T_{2n+1}\le\tau$ and $T_{2n}\ge\tau\,,$ $\forall \ n\in\nat\,.$
Then, for
 $t\in I_{2n+1},$
\begin{equation}\label{Pi1}
\frac {d }{dt}\Vert U(t)\Vert^2\le
B_{2n+1}\Vert U(t)\Vert^2+
B_{2n+1}\Vert U(t_{2n})\Vert^2\,.
\end{equation}
\end{Proposition}

\noindent{\bf Proof:} By differentiating $\Vert U(t)\Vert^2$ we get
\begin{eqnarray*}
\frac {d }{dt}\Vert U(t)\Vert^2
= 2\langle U(t), {\mathcal A}U(t)\rangle + 2
\langle U(t), {\mathcal B}(t) U(t-\tau )\rangle
\,.
\end{eqnarray*}
Then, by using the dissipativness of the operator ${\mathcal A},$
\begin{eqnarray*}
\frac {d }{dt}\Vert U(t)\Vert^2
\le 2
\langle U(t), {\mathcal B}(t) U(t-\tau )\rangle
\,.
\end{eqnarray*}

Hence,
from ${\rm 2}),$
$$ \frac {d }{dt}\Vert U(t)\Vert^2
\le {B_{2n+1}}\Vert U(t)\Vert^2+{B_{2n+1}}\Vert U(t-\tau )\Vert^2\,.
$$
We can now  conclude observing that since $T_{2n+1}\le\tau$ and $T_{2n}\ge\tau\,,$ then for $t\in I_{2n+1}$ it is $t-\tau \in I_{2n}.$
Then, since $\Vert U(t)\Vert$ is decreasing in $I_{2n}\,,$ the estimate in the statement is proved.\qed

The stability result follows.
\begin{Theorem}\label{CP3}
Assume $\mbox{\rm 1),\ 2)},$ $T_{2n+1}\le\tau$ and $T_{2n}\ge\tau\,,$
$\forall\ n\in \nat\,.$
Moreover
assume $T_{2n}> T^*,$ for all $n\in\nat,$ where $T^*$ is the time defined in $(\ref{Tstar}).$
If
\begin{equation}\label{star3}
\sum_{n=0}^\infty
\ln \left [
e^{B_{2n+1}T_{2n+1}}
 (c_n+1-e^{-B_{2n+1} T_{2n+1}}
  )\right ] = -\infty\,,
\end{equation}
then every solution  $U$ of $(\ref{problema})$ satisfies
$\Vert U(t)\Vert \rightarrow 0$ for $t\rightarrow +\infty\,.$
\end{Theorem}

\noindent {\bf Proof.} For $t\in I_{2n+1}=[t_{2n+1}, t_{2n+2}),$
from estimate (\ref{Pi1}) we have
$$
\Vert U(t)\Vert^2\le e^{B_{2n+1}(t-t_{2n+1})}\Big \{ \Vert U(t_{2n+1})\Vert^2 +B_{2n+1} \int_{t_{2n+1}}^t \Vert U(t_{2n})\Vert^2 e^{-B_{2n+1}(s-t_{2n+1})} ds\Big \}.
$$

\noindent
Then we deduce
$$\Vert U(t)\Vert^2
\le e^{B_{2n+1}T_{2n+1} }\Vert U(t_{2n+1})\Vert^2+
e^{B_{2n+1}(t-t_{2n+1})}\Vert U(t_{2n})\Vert^2\Big [ 1-e^{-B_{2n+1}(t-t_{2n+1})}\Big ]\,,
$$
and therefore
$$\Vert U(t)\Vert^2\le e^{B_{2n+1}T_{2n+1} }\Vert U(t_{2n+1})\Vert^2 +e^{B_{2n+1}T_{2n+1} }\Vert U(t_{2n})\Vert^2 -\Vert U(t_{2n})\Vert^2\,,$$
for $t\in I_{2n+1}=[t_{2n+1},t_{2n+2}),\ n\in\nat.$

Now we use the estimate (\ref{stimaobs}) obtaining

$$\Vert U(t_{2n+2})\Vert^2\le e^{B_{2n+1} T_{2n+1}}\left (  c_n +1-e^{-B_{2n+1} T_{2n+1}}
\right )
\Vert U(t_{2n})\Vert^2,\quad n\in\nat\,.$$
Thus,
\begin{equation}\label{Pi2}
\Vert U(t_{2n+2})\Vert \le \Big [
\displaystyle{
\Pi_{k=0}^n}e^{B_{2k+1} T_{2k+1}} (c_k +1 -e^{-B_{2k+1} T_{2k+1}})
\Big ]^{\frac 12} \Vert U_0\Vert \,.\end{equation}
Then the asymptotic stability result follows  if
$$
 \displaystyle{
\Pi_{k=0}^n}e^{B_{2k+1} T_{2k+1}} \left ( c_k +1 -e^{-B_{2k+1} T_{2k+1}}\right )
 \rightarrow 0,\quad \mbox{\rm for}\ n\rightarrow\infty\,,
$$
namely if

$$
 \sum_{n=0}^\infty
\ln \Big [e^{B_{2n+1} T_{2n+1}} ( c_n +1 -e^{-B_{2n+1} T_{2n+1}})
\Big ] \rightarrow -\infty,\quad \mbox{\rm for}\ n\rightarrow\infty\,.\quad\quad
\qed$$

\begin{Remark}{\rm
Observe that, when the odd intervals $I_{2n+1}$ have  length lower or equal than the time delay $\tau\,,$ the assumption (\ref{star3}) is a bit less restrictive than (\ref{general})\,.
Indeed, 
$$e^{B_{2n+1}T_{2n+1}}(c_n+1-e^{-B_{2n+1}T_{2n+1}})
<e^{2b_{2n+1}T_{2n+1}}(c_n+B_{2n+1}T_{2n+1}),\quad \forall n\in \nat\,.$$
}
\end{Remark}

\begin{Remark}{\rm Arguing as in Remark \ref{particular} we can show that
condition (\ref{star3}) is verified, in particular, if (\ref{star1bis}) holds true.}
\end{Remark}

\section{Exponential stability
\label{stexp}}

\hspace{5mm}

\setcounter{equation}{0}
Under additional assumptions on the coefficients $T_n, B_{2n+1}, c_n,$  exponential stability  results hold true.

\begin{Theorem}\label{expTh}
Assume $\mbox{\rm 1), 2)}.$
Moreover, assume
\begin{equation}\label{E1}
T_{2n}=T^0\quad\forall\ n\in\nat ,
\end{equation}
with $T^0\ge \tau$ and $T^0>T^*,$ where  $T^*$ is the constant defined in $(\ref{Tstar}),$
\begin{equation}\label{E2}
T_{2n+1}=\tilde T\quad\forall\ n\in\nat 
\end{equation}
and
\begin{equation}\label{ASS1A}
\sup_{n\in\nat} \ e^{2 B_{2n+1}\tilde T}(c+B_{2n+1}\tilde T)=d<1,
\end{equation}
where $c=Me^{-\mu T^0}.$ Then, there exist two positive constants $C,\alpha$
such that
\begin{equation}\label{expestimateA}
\Vert U(t)\Vert \le C e^{-\alpha t} \Vert U_0\Vert ,\ \ t>0,
\end{equation}
for any solution of problem  $(\ref{problema}).$
\end{Theorem}

\noindent {\bf Proof.} 
Note that, from the definition of the constant $c,$ estimate (\ref{stimaobs}) holds with $c_n=c,\ \forall\ n\in\nat.$ Thus,  from (\ref{ASS1A}) and (\ref{cri4})
we obtain
$$\Vert U(T^0+\tilde T)\Vert \le d ^{\frac 12}\Vert U_0\Vert ,$$
and then, 
$$\Vert U(n(T^0+\tilde T))\Vert\le d^{\frac n2}\Vert U_0\Vert  ,\quad \forall n\in\nat.$$
Therefore, $\Vert U(t)\Vert$ satisfies an exponential estimate like
(\ref{expestimateA}) (see Lemma 1 of \cite{Gugat}).\qed


Concerning the case of  {\em small} delay intervals, namely
$\vert I_{2n+1}\vert \le \tau ,\ \forall n\in \nat\,,$ one can state the following asymptotic stability result.

\begin{Theorem}\label{exp2}
Assume $\mbox{\rm 1), 2)}.$ Moreover assume
$$
T_{2n}=T^0\quad\forall\ n\in\nat ,
$$
with $T^0\ge \tau$ and $T^0>T^*,$ where the time $T^*$ is defined in   $(\ref{Tstar}),$
\begin{equation}\label{E2bis}
T_{2n+1}=\tilde T, \quad \mbox{\rm with}\ \ \tilde T\le\tau \quad\forall\ n\in\nat 
\end{equation}
and
\begin{equation}\label{ASS1Abis}
\sup_{n\in\nat} \ e^{ B_{2n+1}\tilde T}(c+1 -e^{-B_{2n+1}\tilde T})=d<1,
\end{equation}
where $c=Me^{-\mu T^0}.$ Then, there exist two positive constants $C,\alpha$
such that
\begin{equation}\label{expestimateAbis}
\Vert U(t)\Vert \le C e^{-\alpha t} \Vert U_0\Vert ,\ \ t>0,
\end{equation}
for any solution of  $(\ref{problema}).$
\end{Theorem}

\noindent {\bf Proof.} The proof is analogous to the one of Theorem \ref{expTh}\,.\qed

\section{Intermittent anti--damping\label{anti}}

\hspace{5mm}

\setcounter{equation}{0}

With analogous technics we can also deal with an intermittent anti--damping term.
More precisely, let us consider the model 
 \begin{equation}\label{AD}
\left\{
\begin{array}{l}
U_t(t)={\mathcal A}U(t)+{\mathcal B}(t) U(t) \quad t >0,\\
U(0)=U_0,
\end{array}
\right.
\end{equation}
where $\tau$ is the time delay,  the initial datum $U_0$ belongs to ${\mathcal H}$  and, for $t>0,$ ${\mathcal B}(t)$ is a bounded operator from ${\mathcal H}$ such that 
$$\langle {\mathcal B}(t) U, U \rangle \ge 0,\quad\forall \ U\in {\mathcal H}.$$
Thus ${\mathcal B}(t)U(t)$ is  an anti--damping term (cfr. \cite{HMV}). In particular
we consider an intermittent feedback, 
that is
we assume 
that there exists an increasing sequence of positive real numbers $\{t_n\}_n,$ with $t_0=0,$ 
such that
$$
\begin{array}{l}
{ 3)} \quad {\mathcal B}(t)=0\quad  \forall\ t\in I_{2n}=[t_{2n},t_{2n+1}),\\
{4)} \quad  {\mathcal B}(t)={\mathcal D}_{2n+1} 
\quad \forall \  t\in I_{2n+1}=[t_{2n+1},t_{2n+2}).
\end{array}
$$
We denote $D_{2n+1}=\Vert {\mathcal D}_{2n+1} 
\Vert_{{\mathcal L}({\mathcal H})}\,,$ $n\in\nat\,.$

 As before, denote by $T_n$ the length of the interval $I_n,$ that is
$$
T_n=t_{n+1}-t_n,\quad n\in \nat\,.
$$

Note that Proposition \ref{obs}, which gives an observability estimate on the intervals $I_{2n}$ where the anti-damping is not present, still holds true.
Concerning the time intervals $I_{2n+1}$ where the anti-damping acts one can obtain the following estimate.

\begin{Proposition}\label{badAD}
Assume $3)$ and  $4).$ 
For every solution of problem $(\ref{AD}),$ 
$$\frac {d}{dt}\Vert U(t)\Vert^2 \le  2 D_{2n+1}  \Vert U(t)\Vert ^2, \quad t\in I_{2n+1}=[t_{2n+1}, t_{2n+2}], \ \ \forall \ n\in \nat\,.$$
\end{Proposition}

\noindent {\bf Proof.} Being ${\mathcal A}$ dissipative,
the estimate follows immediately from 3).\qed 

From Proposition \ref{badAD} we deduce an asymptotic stability result.

\begin{Theorem}\label{CP1AD}
Assume $\mbox{\rm 3), 4)}.$ Moreover
assume $T_{2n}> T^*,$ for all $n\in\nat,$ where $T^*$ is the time defined in $(\ref{Tstar}).$ If
\begin{equation}\label{generalAD}
\sum_{n=0}^\infty
\ln \left (
e^{2D_{2n+1}T_{2n+1}}
 c_n\right ) = -\infty\,,
\end{equation}
then the
problem  $(\ref{AD})$ is asymptotically stable, that is any solution  $U$ of $(\ref{AD})$ satisfies
$\Vert U(t)\Vert \rightarrow 0$ for $t\rightarrow +\infty\,.$
\end{Theorem}

\noindent {\bf Proof.}
From Proposition \ref{badAD} we have  the estimate
$$
\frac {d}{dt}\Vert U(t)\Vert^2 \le 2  D_{2n+1}  \Vert U(t)\Vert ^2, \quad t\in I_{2n+1}=[t_{2n+1}, t_{2n+2}], \ \ \forall \ n\in \nat\,.$$
This implies
\begin{equation}\label{cri0A}
\Vert U(t_{2n+2})\Vert^2 \le e^{2D_{2n+1}T_{2n+1}}\Vert U(t_{2n+1})\Vert^2 ,
\quad \forall \ n\in\nat.
\end{equation}
Then, from estimate (\ref{stimaobs}) which is always valid of course in the time intervals 
without damping, 
\begin{equation}\label{cri4A}
\Vert U(t_{2n+2})\Vert^2\le e^{2D_{2n+1}T_{2n+1}}
 c_n
\Vert U(t_{2n})\Vert^2 ,
\quad \forall \ n\in\nat\,.
\end{equation}
By repeating this argument  we obtain

\begin{equation}\label{cri5R}
\Vert U(t_{2n+2})\Vert^2 \le\displaystyle{\Pi_{k=0}^n }e^{2D_{2k+1}T_{2k+1}}
 c_k\Vert U_0\Vert^2 ,
\quad \forall \ n\in\nat\,.
\end{equation}
Therefore, asymptotic stability is ensured if
$$\displaystyle{
\Pi_{k=0}^n} e^{2D_{2k+1}T_{2k+1}}
 c_k \longrightarrow 0,\quad \mbox{\rm for}\ n\rightarrow\infty,$$
 or equivalently
 $$\ln \Big (
\Pi_{k=0}^n e^{2D_{2k+1}T_{2k+1}}
 c_k\Big ) \longrightarrow -\infty,\quad \mbox{\rm for}\ n\rightarrow\infty\,.$$
 This concludes.  \qed

\begin{Remark}\label{particularAD}
{\rm
In particular
 (\ref{generalAD}) is verified under the following assumptions:
\begin{equation}\label{star1bisAD}
\sum_{n=0}^\infty D_{2n+1}T_{2n+1}<+\infty\quad \mbox{\rm and}\quad \sum_{n=0}^\infty \ln c_n=-\infty\,.
\end{equation}
}
\end{Remark}

Under additional assumptions on the problem coefficients $T_n, D_{2n+1}, c_n,$ an exponential stability  result holds.

\begin{Theorem}\label{expAD}
Assume $\mbox{\rm 3), 4)}$
and 
\begin{equation}\label{E1AD}
T_{2n}=T^0\quad\forall\ n\in\nat ,
\end{equation}
with $T^0>T^*,$ where the time $T^*$ is defined in $(\ref{Tstar}).$
Assume also that
\begin{equation}\label{E2AD}
T_{2n+1}=\tilde T\quad\forall\ n\in\nat 
\end{equation}
and
\begin{equation}\label{ASS1AAD}
\sup_{n\in\nat} \ e^{2D_{2n+1}\tilde T}c=d<1,
\end{equation}
where, $c=Me^{-\mu T^0}.$  Then, there exist two positive constants $C ,\alpha$
such that
\begin{equation}\label{expestimateAAD}
\Vert U(t)\Vert \le C e^{-\alpha t} \Vert U_0\Vert ,\ \ t>0,
\end{equation}
for any solution of problem  $(\ref{AD}).$
\end{Theorem}

\section{Concrete examples\label
{esempi}}
\hspace{5mm}

\setcounter{equation}{0}
In this section we illustrate some eaxamples falling into the previous abstract setting.

\subsection{Viscoelastic wave type equation \label
{visco}}
\hspace{5mm}

Let $H$ be a real Hilbert space and let $A:{\mathcal D}(A)\rightarrow H$
be a positive self--adjoint  operator with a compact inverse in $H.$ Denote by $V:={\mathcal D}(A^{\frac 1 2})$ the domain of
$A^{\frac 1 2}.$

Let us consider the problem
 \begin{eqnarray}
& &u_{tt}(x,t) +A u (x,t)-\int_0^\infty \mu(s) A u(x, t-s) ds +b(t) u_t(x, t-\tau) =0\quad t>0,\label{1.1}\\
& & u(x,t) = 0 \quad \mbox{\rm on}\quad \partial\Omega\times (0,+\infty ),\label{CB}\\
& &u(x,t)=u_0(x,t)\quad \mbox{\rm in}\quad \Omega\times (-\infty, 0];\label{1.2}
\end{eqnarray}

where  the initial datum $u_0$ belongs to a suitable space, the constant $\tau >0$ is the time delay, and
the  memory kernel $\mu :[0,+\infty)\rightarrow [0,+\infty)$ satisfies

i) $\mu\in C^1 (\RR^+ )\cap L^1(\RR^+ );$

ii) $\mu (0)=\mu_0>0;$

iii) $\int_0^{+\infty} \mu (t) dt=\tilde \mu <1;$

iv) $\mu^{\prime} (t)\le -\delta \mu (t), \quad \mbox{for some}\ \ \delta >0.$

\noindent Moreover, the function $b(\cdot )\in L^\infty_{loc}(0,+\infty )$ is a function which is zero intermittently. That is,
we assume that for all $n\in\nat$ there exists $t_n>0,$ with $t_0=0$ and $t_n<t_{n+1},$
such that
$$
\begin{array}{l}
{ 1_w)} \quad  b(t)=0\quad  \forall\ t\in I_{2n}=[t_{2n},t_{2n+1}),\\
{2_w)} \quad \vert b(t)\vert \le b_{2n+1}\ne 0
\quad \forall \  t\in I_{2n+1}=[t_{2n+1},t_{2n+2}).
\end{array}
$$

Stability result for the above problem were firstly obtained in \cite{pignotti2016}. We want to show that they can also be obtained as particular case of previous abstract setting.  

To this aim, following  Dafermos \cite{Dafermos}, we can  introduce the new variable
\begin{equation}\label{eta}
\eta^t(x,s):=u(x,t)-u(x,t-s).
\end{equation}

Then, problem (\ref{1.1})--(\ref{1.2}) may be rewritten as

\begin{eqnarray}
& &u_{tt}(x,t)= -(1-\tilde \mu)A u (x,t)-
\int_0^\infty \mu (s) A \eta^t(x,s) ds\nonumber\\
& &\hspace{5 cm}
 -b(t)u_t(x,t-\tau )\quad \mbox{\rm in}\ \Omega\times
(0,+\infty),\label{e1d}\\
& & \eta_t^t(x,s)=-\eta^t_s(x,s)+u_t(x,t)\quad \mbox{\rm in}\ \Omega\times
(0,+\infty)\times (0,+\infty ),\label{e2d}\\
& &u (x,t) =0\quad \mbox{\rm on}\ \partial\Omega\times
(0,+\infty),\label{e3d}\\
& &\eta^t (x,s) =0\quad \mbox{\rm in}\ \partial\Omega\times
(0,+\infty), \ t\ge 0,\label{e4d}\\
& &u(x,0)=u_0(x)\quad \mbox{\rm and}\quad u_t(x,0)=u_1(x)\quad \hbox{\rm
in}\ \Omega,\label{e5d}\\
& & \eta^0(x,s)=\eta_0(x,s) \quad \mbox{\rm in}\ \Omega\times
(0,+\infty), \label{e6d}
\end{eqnarray}
where
\begin{equation}\label{datiinizd}
\begin{array}{l}
u_0(x)=u_0(x,0), \quad x\in\Omega,\\
u_1(x)=\frac {\partial u_0}{\partial t}(x,t)\vert_{t=0},\quad x\in\Omega,\\
\eta_0(x,s)=u_0(x,0)-u_0(x,-s),\quad x\in\Omega,\  s\in (0,+\infty).
\end{array}
\end{equation}
Set
$L^2_{\mu}((0, \infty); V)$  the Hilbert space
of $V-$ valued functions on $(0,+\infty),$
endowed with the inner product
$$\langle \varphi, \psi\rangle_{L^2_{\mu}((0, \infty);V)}=
\int_0^\infty \mu (s)\langle A^{1/2} \varphi (s), A^{1/2} \psi (s) \rangle_H ds\,.
$$

\noindent
Let ${\mathcal H}$ be
the Hilbert space

$${\mathcal H}=
V\times H\times L^2_{\mu}((0, \infty); V),$$
equipped
  with the inner product

\begin{equation}\label{innerd}
\begin{array}{l}
\left\langle
\left (
\begin{array}{l}
u\\
v\\
w
\end{array}
\right ),\left (
\begin{array}{l}
\tilde u\\
\tilde v\\
\tilde w
\end{array}
\right )
\right\rangle_{\mathcal H}
:= \displaystyle{
 (1-\tilde\mu )\langle A^{1/2} u, A^{1/2}\tilde u \rangle_H + \langle v, \tilde v \rangle_H +
 \int_0^\infty \mu (s)\langle A^{1/2} w, A^{1/2}\tilde w \rangle_H ds\,. }
\end{array}
\end{equation}
Denoting by $U$ the vector $U=(u, u_t,\eta ),$ the above problem can be rewritten in the form  (\ref{problema}), 
where ${\mathcal B}U=B (u, v, \eta )= (0, bv, 0)$ and
${\mathcal A}$ is defined by
\begin{equation}\label{Operator}
{\mathcal A}\left (
\begin{array}{l}
u\\v\\w
\end{array}
\right )
:=\left (
\begin{array}{l}
v\\
(1-\tilde\mu)A u+\int_0^{\infty}\mu (s) A w(s)ds\\
-w_s+v
\end{array}
\right )\,,
\end{equation}
with domain (cfr. \cite{Pata})

\begin{equation}\label{dominioOpd}
\begin{array}{l}
{\mathcal D}({\mathcal A}):=\left\{
\ (u,v,\eta )^T\in   H^1_0(\Omega)\times H^1_0(\Omega)
\times L^2_{\mu}((0,+\infty);H^1_0(\Omega))\, :\right.\\\medskip
\hspace{3 cm}\left.
(1-\tilde\mu)u+\int_0^\infty \mu (s)\eta (s) ds \in H^2(\Omega)\cap H^1_0(\Omega),\ \eta_s\in  L^2_{\mu}((0,+\infty);H^1_0(\Omega))
\right\}.
\end{array}
\end{equation}
It has been proved in \cite{GRP} that the above system is exponentially stable, namely that the operator ${\mathcal A}$ generates a strongly continuos semigroup satisfying the estimate 
$(\ref{semigroup}),$ for suitable constants. Moreover, it is well-known that, the operator   
${\mathcal A}$  is dissipative. 
Therefore, our previous results apply to this model.

As a concrete example we can consider the wave equation with memory.
More precisely, let $\Omega\subset\RR^n$ be an open bounded domain   with a smooth  boundary
$\partial\Omega.$
Let us consider the initial boundary value  problem

\begin{eqnarray}
& &u_{tt}(x,t) -\Delta u (x,t)+
\int_0^\infty \mu (s)\Delta u(x,t-s) ds\nonumber \\
& &\hspace{5 cm}
 +
b(t) u_t(x,t-\tau )=0\quad \mbox{\rm in}\ \Omega\times
(0,+\infty),\label{1.1d}\\
& &u (x,t) =0\quad \mbox{\rm on}\ \partial\Omega\times
(0,+\infty),\label{1.2d}\\
& &u(x,t)=u_0(x, t)\quad \hbox{\rm
in}\ \Omega\times (-\infty, 0]. \label{1.3d}
\end{eqnarray}

This problem enters in previous form $(\ref{1.1})-(\ref{1.2}),$ if we take
$H=L^2(\Omega)$ and the operator  $A$ defined by
$$A:{\mathcal D}(A)\rightarrow H\,:\,  u\rightarrow -\Delta u,$$
where ${\mathcal D}(A)=H^1_0(\Omega)\cap
H^2(\Omega).$

The operator $A$ is a self--adjoint and positive operator with a compact inverse in $H$
and is such that $V={\mathcal D}(A^{1/2})=H^1_0(\Omega).$

Under the same conditions that before on the memory kernel $\mu (\cdot)$ and on the function $b(\cdot ),$ previous asymptotic/exponential stability results are valid.
The case $b$ constant has been studied in \cite{AlNP2015}.
In particular, we have proved that the exponential stability is preserved, in presence of the delay feedback, if the coefficient $b$ of this one is sufficiently small.
The choice $b$ constant was made only for the sake of clearness. The result in  \cite{AlNP2015} remains true if instead of $b$ constant we consider $b=b(t),$
under a suitable smallness condition on the $L^\infty-$norm of   $b(\cdot )\,.$
On the contrary here we give stability results without restrictions on the  $L^\infty-$norm of   $b(\cdot )\,,$
even if only for on--off $b(\cdot )\,.$

Our results also apply to Petrovsky  system with viscoelastic damping  with Dirichlet and Neumann boundary conditions:

\begin{eqnarray}
& &u_{tt}(x,t) +\Delta^2 u (x,t)-
\int_0^\infty \mu (s)\Delta^2 u(x,t-s) ds\nonumber \\
& &\hspace{5 cm}
 +
b(t) u_t(x,t-\tau )=0\quad \mbox{\rm in}\ \Omega\times
(0,+\infty),\label{1.1P}\\
& &u (x,t) =\frac {\partial u}{\partial\nu }=0\quad \mbox{\rm on}\ \partial\Omega\times
(0,+\infty),\label{1.2P}\\
& &u(x,t)=u_0(x, t)\quad \hbox{\rm
in}\ \Omega\times (-\infty, 0]. \label{1.3P}
\end{eqnarray}

This problem enters into the previous abstract framework, if we take
$H=L^2(\Omega)$ and the operator  $A$ defined by
$$A:{\mathcal D}(A)\rightarrow H\,:\,  u\rightarrow \Delta^2 u,$$
where ${\mathcal D}(A)=H^2_0(\Omega)\cap
H^4(\Omega),$ with
$$H^2_0(\Omega )=\Big\{
v\in H^2(\Omega )\ :\ v=\frac {\partial v}{\partial\nu }=0 \ \ \mbox{\rm on}\ \ \partial\Omega\
\Big\}\,.$$
The operator $A$ is a self--adjoint and positive operator with a compact inverse in $H$
and is such that $V={\mathcal D}(A^{1/2})=H^2_0(\Omega).$

 Therefore, under the same conditions that before on the memory kernel $\mu (\cdot)$ and on the function $b(\cdot ),$ previous asymptotic/exponential stability results are valid.
\bigskip

\subsection{Locally damped wave equation  equation \label
{esempiwave}}
\hspace{5mm}

Here we consider the wave equation with local internal damping
and intermittent  delay feedbck.
More precisely, let $\Omega\subset\RR^n$ be an open bounded domain   with a boundary
$\partial\Omega$ of class $C^2.$  
Denoting by $m$ the standard multiplier 
$m(x)=x-x_0,\ x_0\in\RR^n,$ 
let $\omega_1$ be the intersection of $\Omega$ with an open neighborhood of the
subset of $\partial\Omega$
\be\label{defgamma0}\Gamma_0=\{\, x\in\partial\Omega\, :\ m(x)\cdot \nu (x)>0\,\}.
\ee

Fixed any subset $\omega_2\subseteq\Omega\,,$ 
let us consider the initial boundary value  problem

 \begin{eqnarray}
& &u_{tt}(x,t) -\Delta u (x,t)+a\chi_{\omega_1} u_t(x,t)+b(t)\chi_{\omega_2} u_t(x,t-\tau)=0\ \mbox{\rm in}\ \Omega\times
(0,+\infty),\label{Wo.1}\\
& &u (x,t) =0\quad \mbox{\rm on}\quad\partial\Omega\times
(0,+\infty),\label{Wo.2}\\
& &u(x,0)=u_0(x)\quad \mbox{\rm and}\quad u_t(x,0)=u_1(x)\quad \hbox{\rm
in}\quad\Omega ,\label{Wo.3}
\end{eqnarray}
where $\chi_{\omega_i}$ denotes the characteristic 
function of $\omega_i,$ $i=1,2,$
$a$ is a positive number   
and
$b$ in $L^\infty (0,+\infty)$ 
is an on--off function satysfying $(1_w)$ and $(2_w)$ of subsection \ref{visco}.
The
initia
datum $(u_0, u_1)$ belongs to  $H^1_0(\Omega)\times L^2(\Omega )\,.$

This problem enters into our previous framework, if we take
$H=L^2(\Omega)$ and the operator  $A$ defined by
$$A:{\mathcal D}(A)\rightarrow H\,:\,  u\rightarrow -\Delta u,$$
where ${\mathcal D}(A)=H^1_0(\Omega)\cap
H^2(\Omega).$  

Now, denoting $U=(u, u_t),$ the problem can be restated in the abstract form (\ref{problema}) where ${\mathcal B}U=B (u, v)= (0, b(t)\chi_{\omega_2}v)$ and
${\mathcal A}$ is defined by
\begin{equation}\label{OperatorlocW}
{\mathcal A}\left (
\begin{array}{l}
u\\v
\end{array}
\right )
:=\left (
\begin{array}{l}
v\\
-A u-a\chi_{\omega_1}v 
\end{array}
\right )\,,
\end{equation}
with domain ${\mathcal D}(A)\times L^2 (\Omega)$ in the Hilbert space ${\mathcal H}=H\times H.$

Concerning the the part without delay feedback, namely the locally damped wave equation
\begin{eqnarray} 
& &u_{tt}(x,t) -\Delta u(x,t)+ a\chi_{\omega_1} u_t(x,t)=0\quad \mbox{\rm in}\quad\Omega\times
(0,+\infty),\label{C.1o}\\
& &u (x,t) =0\quad \mbox{\rm on}\quad\partial\Omega\times
(0,+\infty),\label{C.2o}\\
& &u(x,0)=u_0(x)\quad \mbox{\rm and}\quad u_t(x,0)=u_1(x)\quad \hbox{\rm
in}\quad\Omega ,\label{C.3o}
\end{eqnarray}
it is well--known that, under the previous Lions geometric condition on the set $\omega_1$ (or under the more general assumption of control geometric property \cite{BLR}) where the frictional damping is localized,
an exponential stability result 
 holds  (see e.g. \cite{BLR, Komornikbook, KomLoreti, Lag83, LT, Lions, liu, zuazua}). This
 is equivalent to say that the strongly continuous semigroup generated by the operator ${\mathcal A}$ associated to $(\ref{C.1o})-(\ref{C.3o}),$ namely the one defined in (\ref{OperatorlocW}), satisfies (\ref{semigroup}). As well--known, the operator ${\mathcal A}$ is  dissipative.
 Thus previous abstract stability results are valid also for this model. We emphasize the fact  that the set $\omega_2$ may be any subset 
 of $\Omega,$ not necessarily a subset of $\omega_1\,.$
 On the contrary, in previous stability results for damped wave equation and intermittent delay feedback (see e.g. \cite{JDDE14, GenniCri}) the set $\omega_2$ has to be a subset of $\omega_1\,.$
 On the other hand, now the standard (not delayed) frictional damping is always present in time while in the quoted papers it is on--off like the delay feedback and it acts only on the complementary time intervals with respect to this one.

\bigskip

 {\em E-mail address,}
\\
\quad Cristina Pignotti: \quad{\tt \bf
pignotti@univaq.it}

\end{document}